\documentclass[preprint,12pt]{elsarticle}

\usepackage{amssymb,latexsym,amsmath,amsthm}
\usepackage{graphicx}
\usepackage[usenames]{color}
\usepackage[toc,page]{appendix}
\newtheorem{theorem}{Theorem}[section]
\newtheorem{lemma}[theorem]{Lemma}

\newtheorem{definition}[theorem]{Definition}

\newcommand{\be}[1]{\begin{equation}\label{#1}}
\newcommand{\ee}{\end{equation}}
\newcommand{\bdef}[1]{\begin{definition}\label{#1}}
\newcommand{\bthm}[1]{\begin{theorem}\label{#1}}
\newcommand{\ethm}{\end{theorem}}
\newcommand{\blem}[1]{\begin{lemma}\label{#1}}
\newcommand{\elem}{\end{lemma}}
\newcommand{\ba}[1]{\begin{array}{#1}}
\newcommand{\ea}{\end{array}}
\newcommand{\bi}{\begin{itemize}}
\newcommand{\ei}{\end{itemize}}
\newcommand{\bex}{\begin{exercise}}
\newcommand{\eex}{\end{exercise}}
\newcommand{\ben}{\begin{enumerate}}
\newcommand{\een}{\end{enumerate}}
\newcommand{\bas}{\begin{align*}}
\newcommand{\eas}{\end{align*}}
\newcommand{\bbm}{\begin{bmatrix}}
\newcommand{\ebm}{\end{bmatrix}}
\newcommand{\bvm}{\begin{vmatrix}}
\newcommand{\evm}{\end{vmatrix}}
\newcommand{\bsm}{\left[\begin{smallmatrix}}
\newcommand{\esm}{\end{smallmatrix}\right]}
\newcommand{\bs}[1]{\begin{slide}{#1}}
\newcommand{\es}{\end{slide}}
\newcommand{\bpf}{\begin{proof}}
\newcommand{\epf}{\end{proof}}

\newcommand{\vp}{{\boldsymbol{p}}}

\newcommand{\vt}{{\boldsymbol{t}}}

\newcommand{\vT}{{\boldsymbol{T}}}

\newcommand{\RR}{{\mathbb{R}}}

\newcommand{\makehosquare}[3]%
{\dimen0=#1\advance\dimen0 by -#3%
\vrule height#1 width#2 depth0pt \kern-#2%
\vrule height#1 width#1 depth-\dimen0 \kern-#1%
\vrule height#2 width#1 depth0pt \kern-#3%
\vrule height#1 width#3 depth0pt%
}

\usepackage{amssymb}

\begin{document}

\begin{frontmatter}



\title{Trivariate Spline Representations for Computer Aided Design and Additive Manufacturing\footnote{\copyright This manuscript version is made available under the CC-BY-NC-ND 4.0 license https://creativecommons.org/licenses/by-nc-nd/4.0/} \footnote{https://doi.org/10.1016/j.camwa.2018.08.017}}


\author{Tor Dokken, Vibeke Skytt and Oliver Barrowclough}

\address{SINTEF, P.O. Box 124 Blindern, 0314 Oslo, Norway}

\begin{abstract}
Digital representations targeting design and simulation for Additive Manufacturing (AM) are addressed from the perspective of Computer Aided Geometric Design. We discuss the feasibility for multi-material AM for B-rep based CAD, STL, sculptured triangles as well as trimmed and block-structured trivariate locally refined spline representations. The trivariate spline representations support Isogeometric Analysis (IGA), and topology structures supporting these for CAD, IGA and AM are outlined. The ideas of (Truncated) Hierarchical B-splines, T-splines and LR B-splines are outlined and the approaches are compared. An example from the EC H2020 Factories of the Future Research and Innovation Actions CAxMan illustrates both trimmed and block-structured spline representations for IGA and AM.

\end{abstract}


\end{frontmatter}


\section{Introduction}
\label{Sec-intro}
Today, the additive manufacturing industry is enjoying a boom and the potential impact of AM in the coming years is enormous. The direct market for AM is expected to be \$20 billion by 2020 (McKinsey) \cite{McKinsey}. By 2025 the overall economic impact created by AM is expected to be much higher; reaching \$100 billion to \$250 billion if the industrial implementation continues at the current rate. However, one of the bottlenecks in the adoption of AM is the lack of good tools for design and simulation that address AM directly.

AM is a true born child of digitalization that combines aspects of mathematics, material science, computational sciences and process planning. The ISO/ASTM 5290 \cite{ISO-ASTM-52900} standard defines additive manufacturing (AM) as {\it a process of joining material to make parts from 3D model data, usually layer upon layer, as opposed to subtractive manufacturing and formative manufacturing technologies}. 

The mathematical and computational communities have until recent years paid little attention to AM, and consequently AM technology and research have been addressed mainly from the perspectives of manufacturing and material research. A consequence of this is that the mathematical approaches of Computer Aided Geometric Design (CAGD) from the 1980s still dominate AM, as described in Section \ref{Sec-AM-standards}. In Section \ref{sec-B-spline-basics} we recap univariate B-splines space and tensor product B-splines and use this in Section \ref{sec-tensor-collection} to address spline spaces spanned by collections of tensor product B-splines. Then (Truncated) Hierarchical B-splines, T-splines and LR B-splines are compared in Section \ref{sec-comparing}. In Section \ref{Sec-tri-var-cad} we consider how these novel spline representations can be used in ISO 10303 STEP. The applicability of IGA for analysis based design is discussed in Section \ref{Sec-ana-bas-des}, and an example of block structured and trimmed IGA provided. How to detrim these models for quadrature is addressed in Section \ref{Sec:quadrature}. In Section \ref{Sec:challenges} we address why the approaches outlined in the other sections are not sufficient for the representation of lattice structures.

\section{Standards and object representation for AM}\label{Sec-AM-standards}
The dominant object representations in AM today focus on the needs of manufacturing. Thus, an approximation of the smooth geometric shapes arising in CAD is acceptable if the manufacturing tolerances of the AM-process are respected. In the early days of AM, the STL-format emerged, see Section \ref{Sec-STL}. STL is based on triangulations and targets single material AM. STL is an intermediate and simplified step between B-rep CAD-models (addressed in Section \ref{Sec-B-rep}) and the additive manufacturing process. 

In 2015 the ISO/ASTM 529 Additive Manufacturing File Format (AMF) was introduced providing more accurate shape representation and possibilities of multi-material printing, see Section \ref{Sec-ISO-ASTM}. However, neither STL nor AMF are well aligned with the representation formats of CAD and Finite Element Analysis (FEA). This lack of interoperability makes it very hard to transfer modifications of an object done during the AM process planning back to CAD and analysis-based design.

\subsection{B-rep based CAD and STEP}\label{Sec-B-rep}
When the first 3D printing company, 3D Systems Inc, introduced their Stereolithography technology (SLA) in 1987 \cite{Wohlers-2014}, solid B-rep CAD was emerging \cite{APS-CAD-CAM}. B-rep CAD represents a solid object by its limiting surfaces as it is assumed that the material of the object is uniform. In B-rep CAD, the surface types used are elementary surfaces (e.g., planes, spheres, cones, cylinders, tori) and Non-Uniform Rational B-splines (NURBS). These matched well the subtractive and formative manufacturing process used by industry in the 1980s. 

Data exchange of B-rep CAD models is today mainly by ISO 10303 {\it Automation systems and integration - Product data representation and exchange} (STEP) or through vendor proprietary formats. The development of STEP started in the 1980s, with the intent to make STEP a successor of standards such as IGES, SET and VDA-FS. In 1994/95 ISO published the initial release of STEP as an international standard, thus ending the first phase of the STEP development. STEP is today a widely used standard for exchange of CAD-models and is under continuous development. In 2018 additions addressing trivariate spline representations are expected to be published. These additions can potentially be used as resources for representing object models addressing graded, multi and anisotropic material.

\subsection{STL}\label{Sec-STL}
The Stereolithography technology (SLA) targets objects composed of a single material. Thus B-rep based CAD-models in principle matched the needs of SLA.
In SLA, as well as many other additive processes, each layer is a plane. Consequently, the 3D model data has to be repeatedly sliced by planes. It is algorithmically simpler to slice 3D model data represented by triangulations than the 3D model data using many smooth shape representations. The algorithmic complexity of the slicing is reduced by first tessellating the 3D model data, and then slicing the resulting triangulated model. To represent tessellated models 3D Systems Inc in 1987 introduced the STL-format {\it Standard Triangle Language}, which still today is dominant for AM shape representation. To print objects with high quality, the accuracy of the tessellated model must be adapted to the accuracy of the AM technology used. In addition, the accuracy depends on the scale the object is printed at. An STL file consists only of a list (or soup) of triangles. Each triangle (facet) is represented by its three corners, and an optional normal assigned to the facet. The normal is frequently set to (0,0,0), in which case it is assumed that the direction of the normal can be inferred from the orientation of the vertices of the triangle. There is no guarantee in STL that the collection of triangles represents the closed surface of a volume. Though other tessellation formats such as OBJ, PLY and OFF partially solve this problem by including topology information, other issues with these tessellation formats exist.
Although STL works well for moderately complex objects with single material AM-processes, for objects exhibiting complexity over multiple scales, STL-files are well known to be too voluminous for efficient use. 
It is worth noting that AM is particularly well suited to fabricating models with multiscale complexity, so this issue represents a true bottleneck in the AM pipeline today. 

\subsection{ISO/ASTM 529 - AMF}\label{Sec-ISO-ASTM}
The AM community is aware that more advanced representations are needed, e.g., for supporting multi-material processes and complex lattice structures. For getting a more compact representation of sculptured shapes AMF has introduced a curved triangle patch defined by vertices with optional normals. Sculptured triangles are to be recursively subdivided into four triangles to generate a set of flat triangles to reach the desired resolution. The sculpted triangles are used for describing the surface of the volumes and each volume can be associated with a material ID. An object can consist of a set of volumes each with possible different material IDs. AMF gives requirements that ensure that the geometry represented is well defined. This will drastically improve the quality of model exchange compared to STL. However, as of  2018 the uptake of AMF is slow.

\section{Recap of some B-splines basics}\label{sec-B-spline-basics}
As we will use properties of univariate splines spaces and tensor product B-splines we provide some basics on these topics respectively in Section \ref{sec-univariate-b-splines} and \ref{sec-tensor-product-B-spline}.
\subsection{Univariate splines}\label{sec-univariate-b-splines}Spline representations are popular for curves as they can represent complicated shapes by a sequence of polynomial pieces. They are easy to evaluate meaning that they are suitable for interactive visualization and can thereby support design processes. The continuity between adjacent pieces can be set as required and can be varied according to needs. For many purposes splines of degree three are used. However, the polynomial degree can be chosen according to what is needed for the problem addressed. 

The B-spline basis, $B_{i,p}(t), i=1,\ldots,N$, is very efficient for representing piecewise polynomial curves. It is more accurate than alternative representations, e.g., the representation of each polynomial segment in the power basis. B-splines are defined by a non-decreasing sequence of real numbers denoted knots ${\bf t}=\{t_1,...,t_{N+p+1}\}$. Here $p \ge 0$ is the polynomial degree, $N$ is the dimension of the spline space, and $t_i<t_{i+p+1}, \forall i$. The value of the B-spline is calculated by the recursion relation
\be{eq:B-spline}
\begin{aligned}
B_{i,0}(x)&:=\begin{cases}
1,& t_i\le x<t_{i+1},\\ 0,& \text{otherwise},
\end{cases}\\
B_{i,p}(x)&:=\frac{(x-t_i)}{(t_{i+p}-t_i)}B_{i,p-1}(x) + \frac{(t_{i+p+1}-x)}{(t_{i+p+1}-t_{i+1})}B_{i+1,p-1}(x).
\end{aligned}
\ee

The continuity at each unique knot value is $p-m$ where $m$ is the multiplicity, i.e., the number of times the knot value is repeated. If we want to refine a polynomial curve $f(t)=\sum_{i=1}^N c_i B_{i,p}(t)$ by inserting new knots in the knot sequence ${\bf t}$ then the coefficients for representing $f(t)$ in the refined basis are efficiently calculated by combining weights calculated by the {\it Oslo-algorithm} \cite{Oslo-Algorithm} and the original coefficients of $f(t)$. 

\subsection{Tensor product B-spline }\label{sec-tensor-product-B-spline}
Definition 1.2 in \cite{Dokken-2013} addresses tensor product B-splines. The focus is in single tensor product B-splines and is thus concerned only with the knots in the support of the univariate B-splines multiplied, not the complete knot vector spanning univariate spline spaces.
\bdef{def:Tensor-B-splines}{\bf Tensor product B-splines}.
Let $d$ be a positive integer, suppose $\vp=(p_1,\ldots,p_d)$ has nonnegative integer components, and let $\vt_k:=(t_{k,1},\ldots$, $t_{k,p_k+2})$ be nondecreasing sequences (of knots) $k=1,\ldots,d$. We define a {\bf tensor-product B-spline} $B[\vT]=B[\vt_1,\ldots,\vt_d]:\RR^d\to\RR$ from univariate B-splines $B[\vt_k]$ by
$$B[\vt_1,\ldots,\vt_d](x_1,\ldots,x_d):=\prod_{k=1}^d B[\vt_k](x_k).$$
The {\bf support} of $B$ is given by the cartesian product
\be{eq:tsupp}
\text{\em supp}(B):=[t_{1,1},t_{1,p_1+2}]\times\cdots\times[t_{d,1},t_{d,p_d+2}].
\ee
\end{definition}

\section{Spline spaces spanned by collections of tensor product B-splines}\label{sec-tensor-collection}
A collection of tensor product B-splines will span a spline space. B-spline surfaces in CAD are spanned by a special class of collections of tensor product B-splines. The collection is generated by a tensor product of two univariate spline spaces. Consequently, both B-splines and control points are structured in a regular grid. The regular grid structure was essential when B-splines were introduced in CAD in the 1980s. The tensor product structure of the spline spaces allowed the implementation of very efficient algorithms for interpolation, and evaluation. However, this efficiency comes at the cost of large data increases as the size of the problems increases (both in extent and dimension), as will be discussed later.

Spline spaces that are a tensor product of two univariate spline spaces are central for surface representation in CAD. This relates both to tensor product B-spline surfaces and Non-Uniform Rational B-spline surfaces (NURBS). However, as explained in Section \ref{sec-lack-refinement} such spline spaces lack local refinement. The lack of local refinement is even more severe in the trivariate than in the bivariate case. In AM, trivariate representations are needed when modelling variable material properties and interior structures. 
Due to the lack of local refinement of these representations, there is a demand for more flexible collections of multivariate tensor product B-splines to describe volumes. 

Industrial uptake for an augmented B-spline technology based on more general collections of tensor product B-splines is expected to depend on its compatibility with tensor product B-spline surfaces in CAD. The choice of an augmented spline technology is also feasible to address from the computational perspective. Since the 1960s we have had a doubling in the number of components per integrated circuit every 12 to 18 months (Moore's law). This has in practice given a similar increase in computational power. This computational power allows augmented spline technologies to be explored.

With the above in mind, it seems natural to impose the following requirements on refinement processes for generating more general collections of tensor product B-splines:
\begin{itemize}
\item The starting point is a collection of tensor product B-splines generated by the tensor product of univariate B-spline spaces. These spline spaces are known to span the full polynomial space over each element (polynomial segment).
\item The refinement process creates a nested sequence of refinement spaces. The nesting of the spline spaces will ensure that the full polynomial space is spanned also over each element in the refinement spaces.
\end{itemize}
We will discuss three approaches to refinements following the principles above, namely, Hierarchical B-splines, T-splines and LR B-splines. These will respectively be addressed in Sections \ref{sec-hier}, \ref{sec-T-splines} and \ref{sec-LR}. It should be noted that all the possible spline spaces that can be spanned over box-partitions or T-meshes, cannot be represented by collections of tensor product B-splines. The restriction to spline spaces that can be spanned by a collection of B-splines is computationally feasible and builds on the B-spline technology of state-of-the-art CAD. Alternative approaches to locally refinable splines can be found in \cite{local-survey-2015}. 

\begin{figure}[t]
\begin{center}
\includegraphics[width=13.5 cm]{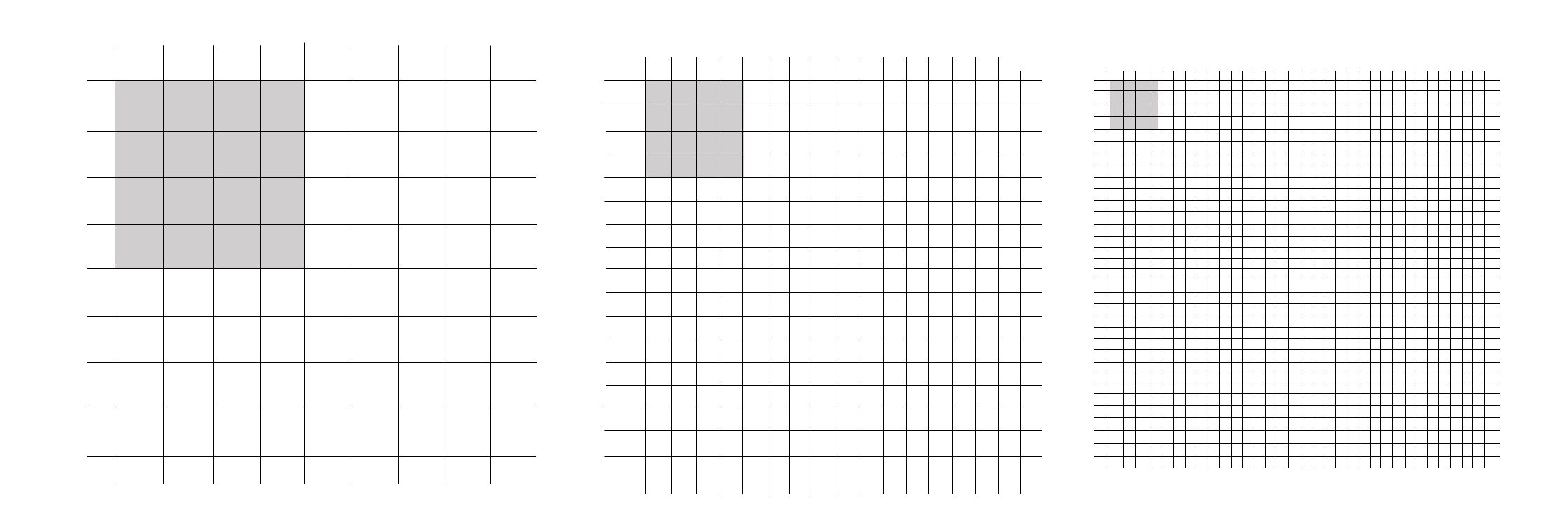}
\caption{To the {\it left} the knotlines of a bi-degree (3,3) uniform tensor product B-splines space at level $l'$ with the support of a sample B-spline. In the {\it middle} knotlines and a sample B-splines at level $l'+1$, and to the 
{\it right} knotlines and a sample B-splines at level $l'+2$.} 
\label{fig:Three_hierarcical_levels}
\end{center}
\end{figure}
\begin{figure}[h]
\begin{center}
\includegraphics[width=13.5 cm]{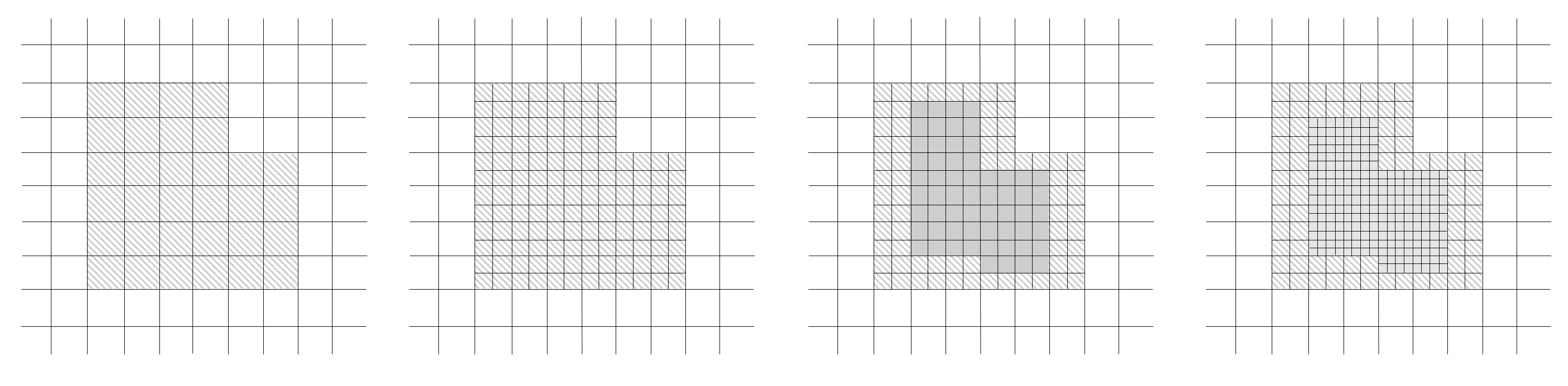}
\caption{To the {\it left} the first mesh from Fig. \ref{fig:Three_hierarcical_levels} with a region that we want to refine with tensor product B-splines from the middle mesh in Fig. \ref{fig:Three_hierarcical_levels}. The second mesh from the left shows the mesh resulting from this refinement. Then in the third mesh 
a region is marked where we want another level of refinement. The mesh to the right shows the final hierarchical mesh.
}\label{fig:Hierarchal-refinement}
\end{center}
\end{figure}

\subsection{Spline spaces spanned by a tensor product of univariate B-spline spaces lack local refinement }\label{sec-lack-refinement}
Traditionally sculptured surfaces in CAD are smooth with limited local variation and not too many degrees of freedom. Consequently, defining a bivariate B-spline basis as a tensor product of two univariate B-spline bases gives an efficient representation. Tensor product B-spline surfaces have been eagerly adopted as a suitable spline representation for sculptured surfaces in CAD. If an extra degree of freedom is needed in the first parameter direction of a tensor product B-spline surface, an extra knot $s$ is inserted. If the surface originally had $N_1 \times N_2$ control points, then the resulting refined surface will have $ N_2$ additional control points. If the extra degree of freedom is needed all along the knot line corresponding to $s$ then the growth by $N_2$ control points is very reasonable. However, if the extra degrees of freedom are only need very locally then the growth is unacceptable. 

Moving to $\mathbb{R}^3$ the issue becomes even more apparent and more server. Let us assume we have a volume spanned by a spline space that is the tensor product of three univariate B-spline spaces respectively of dimension $N_1$, $N_2$ and $N_3$. The collection will contain $N_1 \times N_2 \times N_3$ tensor product B-splines and control points. If we insert an extra knot in the first univariate B-spline space, the number of tensor product B-splines and control points will grow by $N_2 \times N_3$. While the tensor product B-splines can be represented efficiently due to the tensor product of univariate splines space, the coefficients will all have be to be represented, and the corresponding increase in degrees of freedom adds computational complexity, e.g., when solving matrix equations. This lack of local refinement hinders spline spaces that are tensor products of univariate spaces to be used in AM.

\subsection{Hierarchical B-splines}\label{sec-hier}
Hierarchical B-splines (HB) were introduced in 1988 \cite{Forsey:88}. They are based on a dyadic sequence of grids determined by scaled lattices 
$(\frac{k_1}{2^l}, \frac{k_2}{2^l},\ldots,\frac{k_d}{2^l})$. On each of the dyadic grids a tensor product B-spline space with uniform knots is defined as illustrated in Fig. \ref{fig:Three_hierarcical_levels}. 

The refinement is done level by level by removing tensor product B-splines on the coarser level and adding B-splines at the finer level in such a way that linear independence is ensured \cite{Kraft:98}, and the the full polynomial space is spanned over each polynomial element. For an example of refined meshes see Fig. \ref{fig:Hierarchal-refinement}. A partition of unity basis was provided in 2012 when Truncated Hierarchical B-splines (THB) were introduced \cite {Giannelli-2012}. THB-splines provide a basis that is a partition of unity by subtracting scaled tensor product B-splines on the finer level from the tensor product B-splines at the rougher level. It has been observed that some truncated B-splines risk ending up with a support split in two disjoint parts. This problem was addressed by imposing restrictions on allowed refinements \cite{split-domain-THB}. The approach of HB and THB is easily extended to higher dimensions, and any polynomial degree.

\begin{figure}[t]
\begin{center}
\includegraphics[width=13.5 cm]{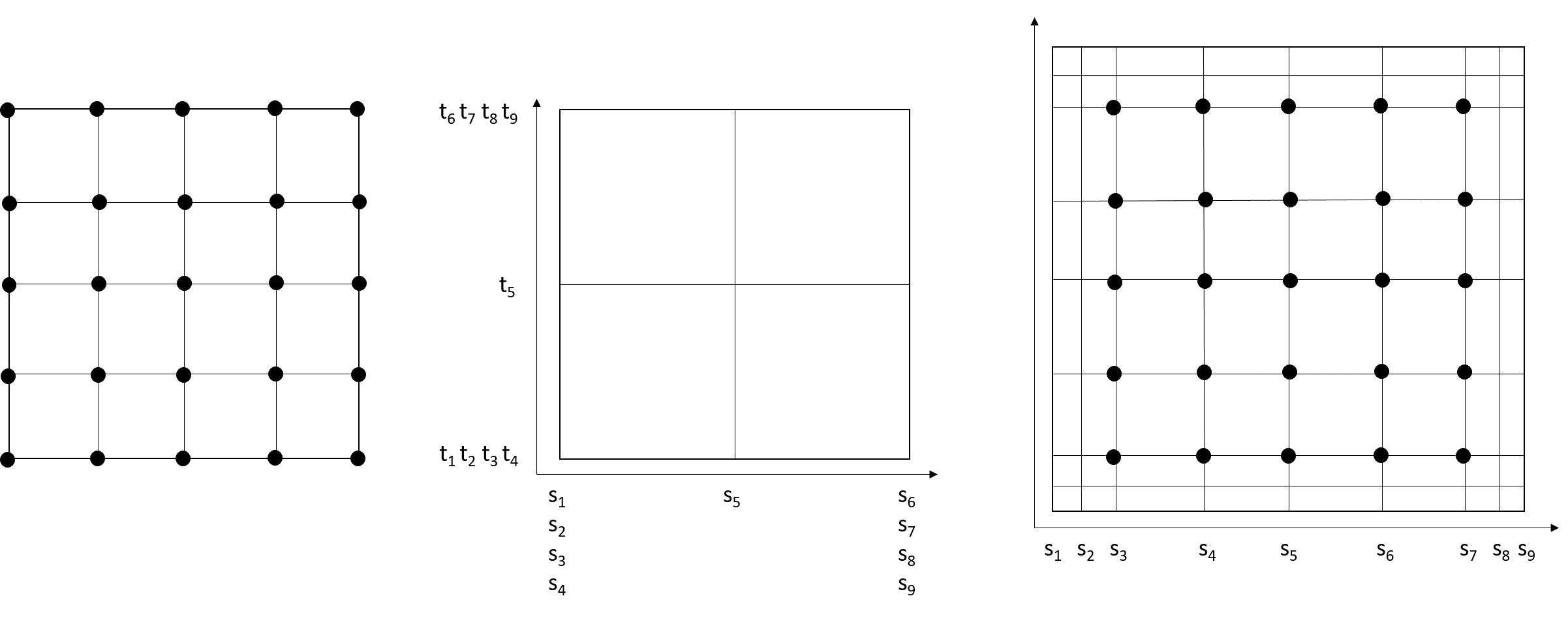}
\caption{Different visualizations of the structure of a bi-degree (3,3) B-spline surface. To the {\it left} we show the control points of the B-spline surface, and in the {\it middle} the piecewise polynomial structure and knotlines. To the {\it right} the control points and knotlines are combined into one graphical representation with no explicit information on knot-multiplicity. We will use the graphical representation to the {\it right} in Fig. \ref{fig:Finding_knots_in_T-mesh} and Fig. \ref{fig:T-mesh_refinement_and_B-spline} to illustrate the idea behind local refinement of T-splines. Note that the pair of middle knots of each tensor product B-spline is anchored at a control point.
}\label{fig:Traditonal_and_combined_B-spline}
\end{center}
\end{figure}

\begin{figure}[h]
\begin{center}
\includegraphics[width=9 cm]{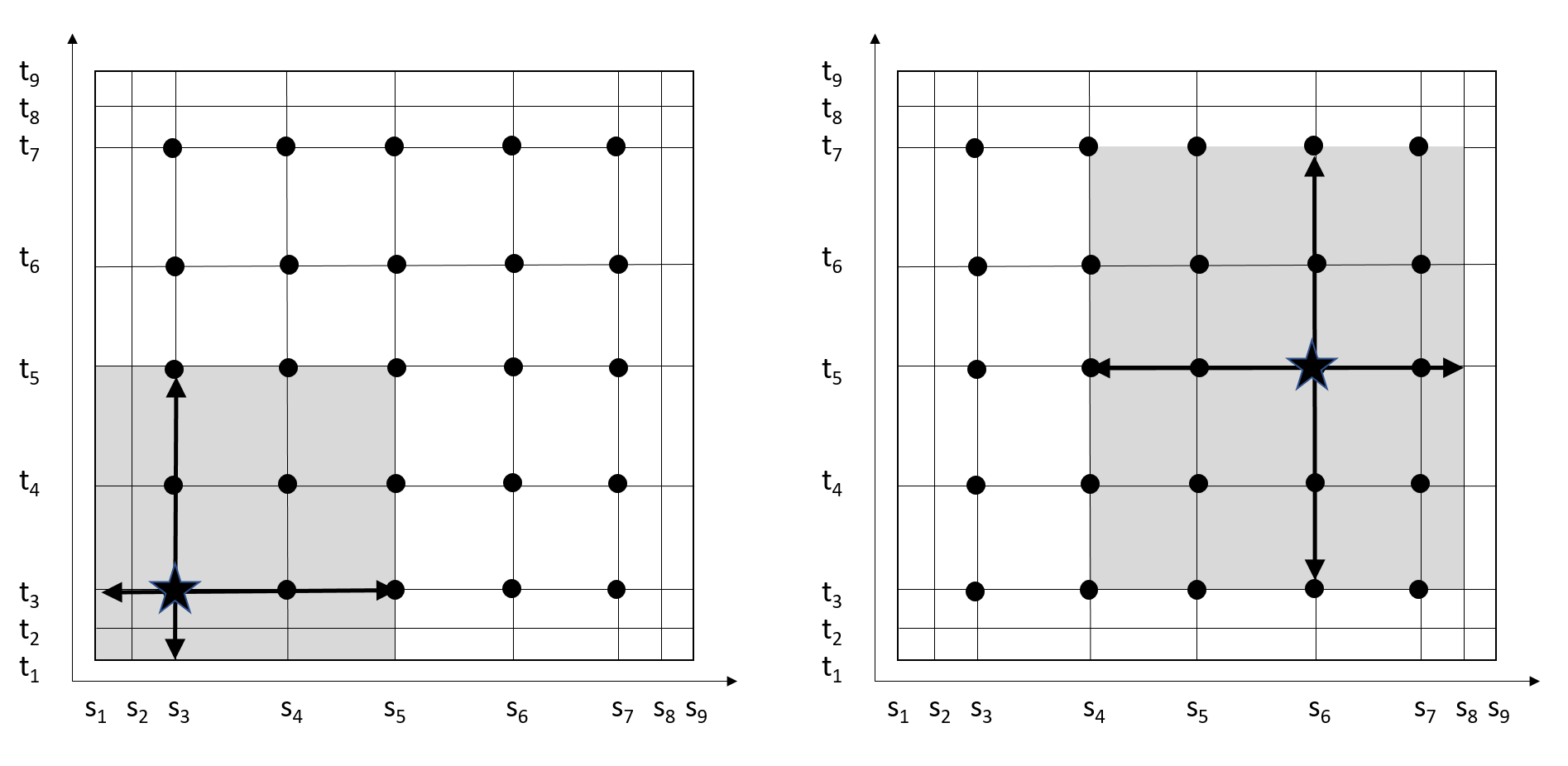}
\caption{The figure to the {\it left} shows that we can identify the knots of the B-spline anchored at $(s_3,t_3)$ by traversing the T-mesh two lines to the left, two lines to the right, two lines down and two lines up. To the {\it right} we show the same for the B-spline anchored at $(s_6,t_5)$.
}\label{fig:Finding_knots_in_T-mesh}
\end{center}
\end{figure}

\begin{figure}[h]
\begin{center}
\includegraphics[width=13.5 cm]{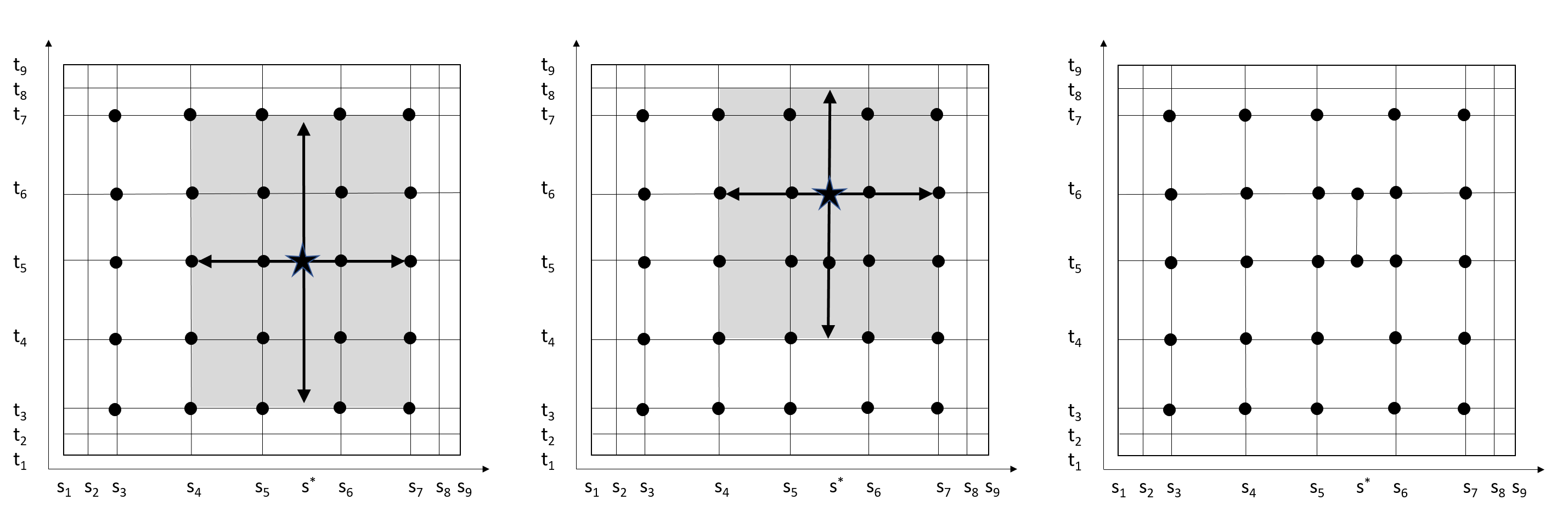}
\caption{To the {\it left} we show that we can add a new control point at $(s^*,t_5)$ and then traverse the mesh to find the knots of the new tensor product B-spline. Note that the knot vector in the first parameter direction of the B-splines anchored at $(s_4,t_5)$, $(s_5,t_5)$, $(s_6,t_5)$ and $(s_7,t_5)$ have to be updated to take $s^*$ into account. These tensor product B-splines all have anchor points with $t_5$ in the second knot direction, and $s^*$ is within their interval of knots in the first knot direction. To the {\it right} we connect the two new control points with a line and make two T-joints, following the second rule of standard T-splines. 
}\label{fig:T-mesh_refinement_and_B-spline}
\end{center}
\end{figure}

\begin{figure}[h]
\begin{center}
\includegraphics[width=10 cm]{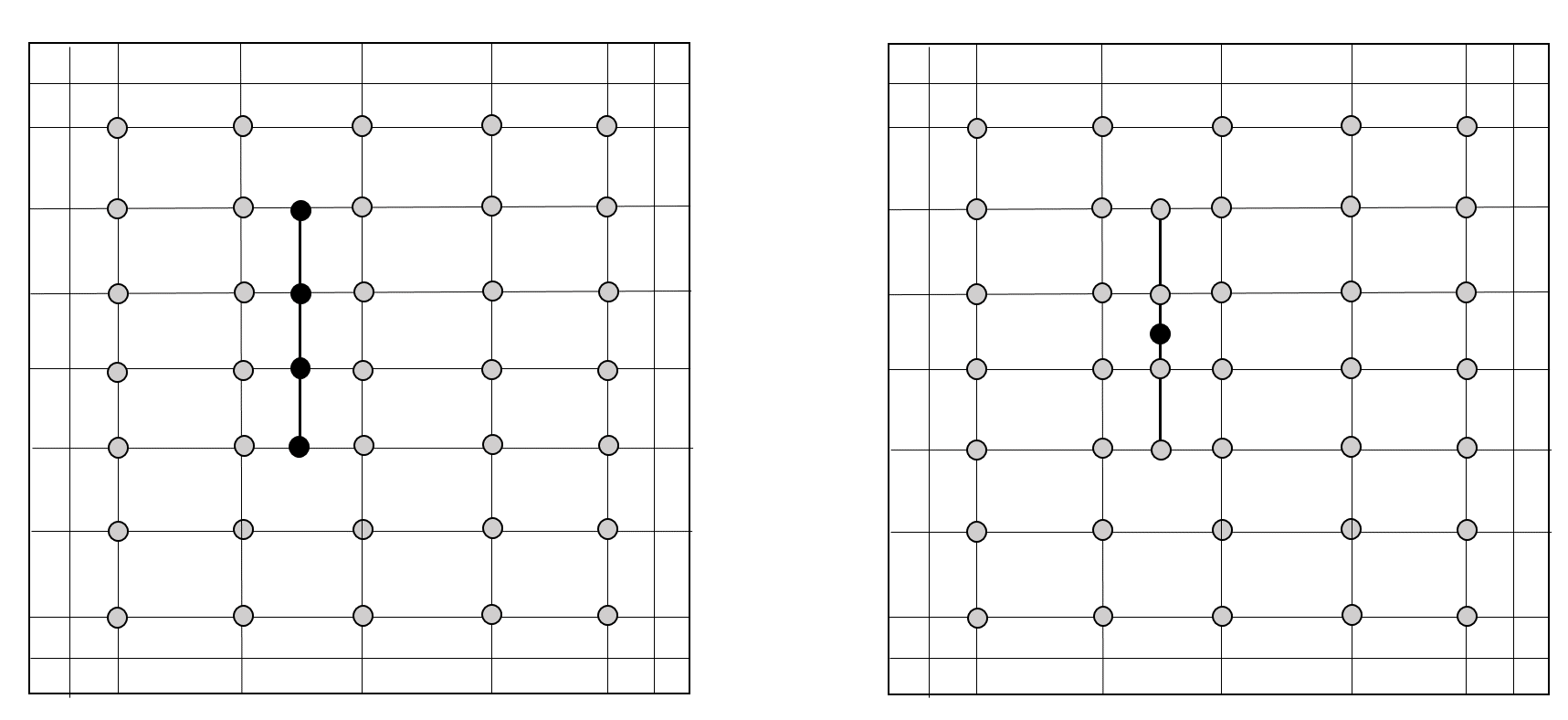}
\caption{To the {\it left} we show that for standard T-splines four adjacent control points in one parameter direction have to be inserted before refinement is allowed in the other parameter direction. Note that the control points all must be on the same constant parameter line. To the {\it right} we show that a control point is inserted in the middle line segment of the line segments created by the insertion of the four control points inserted in the mesh to the {\it left}. 
}\label{fig:Standard_T-splines}
\end{center}
\end{figure}

\begin{figure}[h]
\begin{center}
\includegraphics[width=12.0 cm]{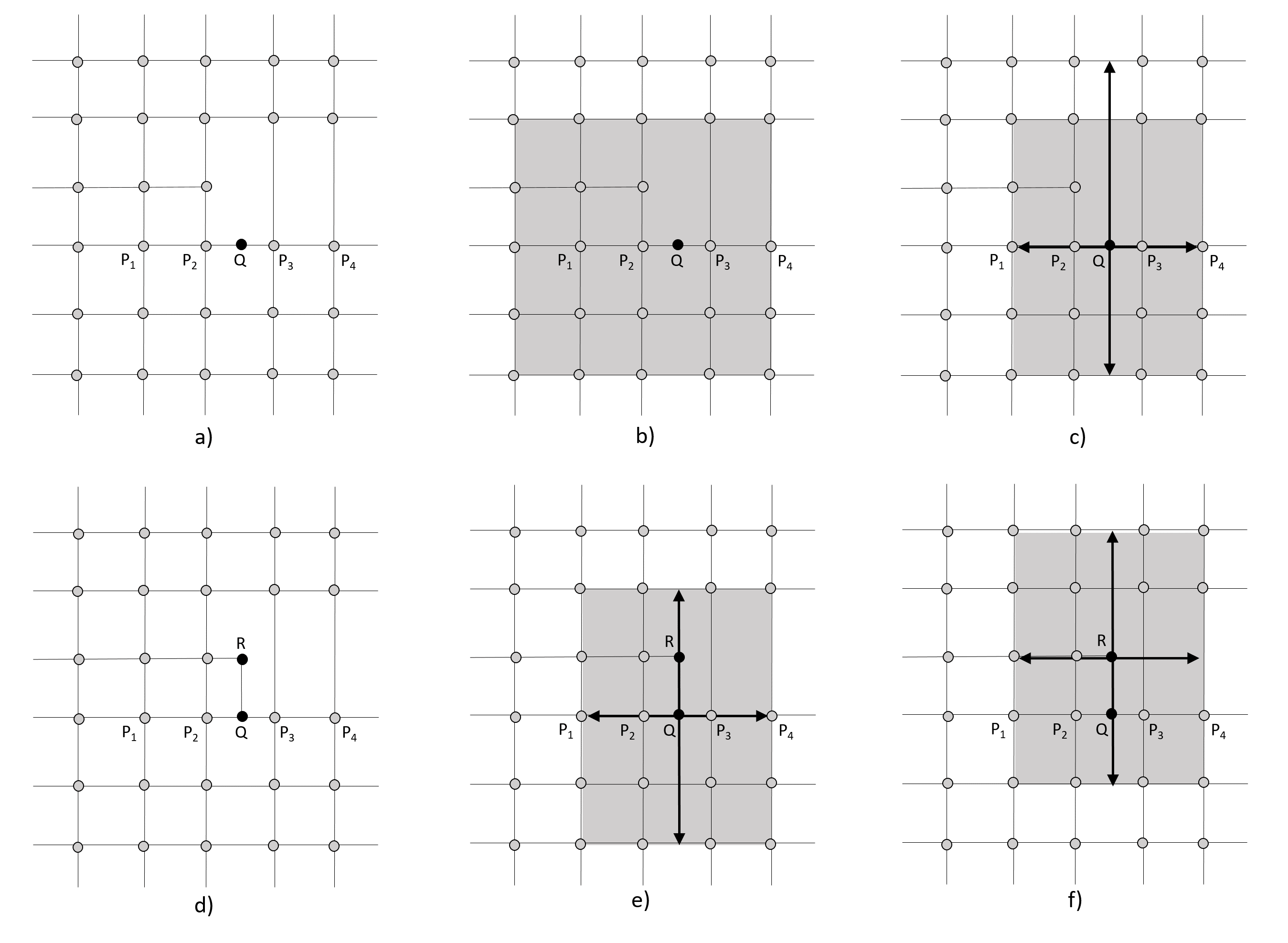}
\caption{In a) we insert a new control point $Q$ in the T-mesh from Figure 9 in \cite{Sederberg:04}. In b) a B-spline to be refined with the first parameter value of $Q$ is shown. The refinement creates a B-spline in c) with centre knots at $Q$. This cannot be found by traversing the mesh. To make this B-spline valid (and to ensure nested spline spaces) a new control point $R$ is added in d), and connected by line segments to existing control points. In e) we see that the B-spline in c) can be generated from the updated mesh. In f) we see an additional B-spline anchored in $R$.}\label{fig:Semi-standard_T-spline}
\end{center}
\end{figure}

\subsection{T-splines}\label{sec-T-splines}

When T-splines were introduced in 2003 the idea was to provide the designer with new interactive tools for local refinement of B-splines surfaces. Isogeometric Analysis introduced in 2005 \cite{Hughes:05} replaced the traditional shape functions of finite elements with B-splines that cross element boarders. It was soon evident that spline spaces made by tensor products of univariate spline spaces were at risk of growing too large for efficient use in IGA. Consequently T-splines gained interest from the IGA community.

T-splines denote a class of locally refined splines, in literature most often presented as bi-degree $(3,3)$, that are refined by successively inserting control points and edges in a so-called {\it T-mesh}. The starting point for T-spline refinement is a bi-degree $(3,3)$ tensor product B-spline surface with control points and knot values organized in a visual rectangular mesh as shown in Fig. \ref{fig:Traditonal_and_combined_B-spline}. In Fig. \ref{fig:Finding_knots_in_T-mesh} we show how this graphical representation can be used for finding the knot values of the tensor product B-spline anchored to each control point. Then in Fig. \ref{fig:T-mesh_refinement_and_B-spline} we insert two new control points and connect these with a line making two T-joints. This illustrates how the rule from Fig. \ref{fig:Finding_knots_in_T-mesh} can be used for finding the knot values of the tensor product B-splines anchored to each of the two new control points. At each of these control points the new line ends at an existing line forming a T-shape. The name T-mesh comes from the T-joints created during such refinement. However, as Fig. \ref{fig:Semi-standard_T-spline} shows, for standard T-splines also L-shapes can occur.

The above procedure is algorithmic and there is no guarantee that nested spline spaces result from the refinement. To ensure that the spline spaces generated are nested additional rules are imposed, defining subclasses of T-splines. In Subsection \ref{sec-standard-T-splines} we will address {\it Standard T-splines} and in Subsection \ref{sec-semi-standard-T-splines} we will address {\it Semi-standard B-splines}. Both these classes create nested spline spaces, and a collection of B-splines that form a (scaled) partition of unity. T-splines that do not form a scaled partition of unity are denoted {\it Non-standard T-splines}. For such T-splines partition of unity is achieved by rational scaling. 

Linear dependence issues detected in \cite{buffa:2010} were avoided by introducing further restrictions on allowed refinements that introduced the subclass of T-splines denoted Analysis Suitable T-splines (AST) \cite{AST:12}. 

T-splines can be directly extended to any odd polynomial degrees, but for even degrees there is no natural middle knot of a B-spline. To rectify this a dual grid is introduced for the anchoring of the tensor product B-splines and control points. The T-spline approach is also used for creating collections of trivariate tensor product B-splines \cite{Solid-T-spline}. Truncation of T-splines has also been proposed \cite{Wei:17}.

\subsubsection{Standard T-splines}\label{sec-standard-T-splines}
The first T-splines introduced were {\it Standard T-splines} \cite{Sederberg:03}. The tensor product B-splines created by standard T-splines form an unscaled partition of unity, $\sum_{\bf i} B_{\bf i}(s,t)=1$. The refinement follows three rules, where the first relates to consistency of knot values, the second relates to when to connect control points (as we did in Fig. \ref{fig:T-mesh_refinement_and_B-spline}). The third restricts when a refinement i allowed: All existing tensor product B-splines that are to be refined must to have identical knot vectors in the other parameter direction, see Fig. \ref{fig:Standard_T-splines}. This can also be regarded as requiring a local tensor product structure on the tensor product B-splines to be updated after insertion of a new control point. For full details consult \cite{Sederberg:03}.

\subsubsection{Semi-standard T-splines}\label{sec-semi-standard-T-splines}

When performing repeated refinement such as shown in Fig. \ref{fig:T-mesh_refinement_and_B-spline} there is no guarantee that the resulting spline space includes prior spline spaces in the refinement sequence. Let $Q$, as in Fig. \ref{fig:Semi-standard_T-spline}.a, be a new control point with knot coordinate $(\sigma,\tau)$ and let $l_\tau$ be the constant knot line on which $Q$ lies. When refining, $Q$ is used for two purposes:
\begin{enumerate}
\item $Q$ is used as an anchor point for a new tensor product B-spline with knots picked from the T-mesh.
\item The knot value $\tau$ is used for refining all B-splines with anchor points on $l_\tau$ and support in the first parameter direction containing $\tau$. 
\end{enumerate}
There is no guarantee that the collection of tensor product B-splines resulting from the above two steps match the T-mesh. (Illustrated in Fig. \ref{fig:Semi-standard_T-spline}.c.) For the B-spline in Fig. \ref{fig:Semi-standard_T-spline}.c to be represented in the T-mesh a new control point $R$ has to be added as shown in Fig. \ref{fig:Semi-standard_T-spline}.d. Now the new spline space is a refinement of the spline space from which we started.

A refinement process such as the above ensures that we generate nested spline spaces and generate {\it Semi-standard T-splines}. A formal algorithm for this process is described in \cite{Sederberg:04}. Below we present a condensed version of the algorithm:
\begin{enumerate}
\item Insert new control points into the T-mesh.
\item Make a collection of B-splines by refining existing tensor product B-splines using $Q$, and creating the tensor product B-splines anchored to the new control points.
\item If any tensor product B-spline misses a knot dictated by mesh traversal, update the tensor product B-spline by knot insertion.
\item If a tensor product B-spline has a knot that is not dictated by mesh traversal, add an appropriate control point in the T-mesh.
\item Repeat 3. and 4. until all issues are solved.
\end{enumerate}

\begin{figure}[h]
\begin{center}
\includegraphics[width=13.5 cm]{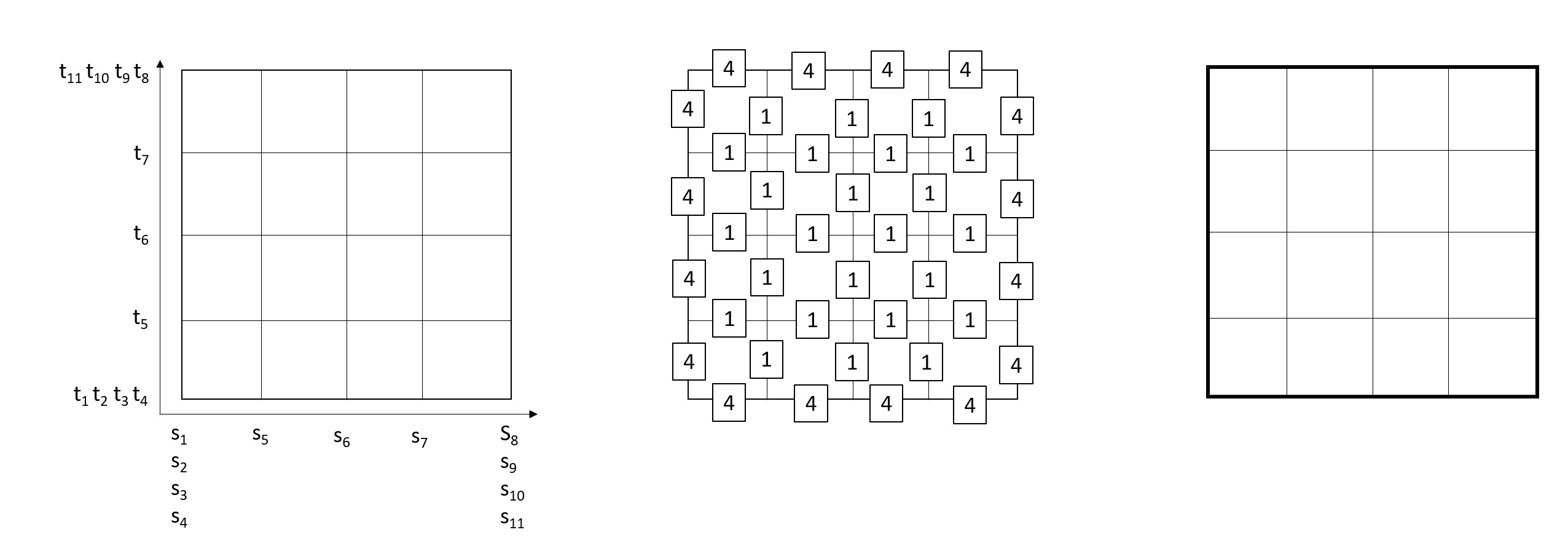}
\caption{To the {\it left} we have the knotlines and knots of a bi-degree (3,3) tensor product basis. In the {\it middle} the knot lines are represented as mesh rectangles with multiplicities assigned. To the {\it right} we represent the meshrectangles of multiplicity 4 with thick lines, and the meshrectangles of multiplicity 1 with thin lines. We will use this convention for meshrectangle multiplicities 1 and 4 in the illustrations to follow. 
}\label{fig:B-spline_to_box-partition}
\end{center}
\end{figure}

\begin{figure}[h]
\begin{center}
\includegraphics[width=13.5 cm]{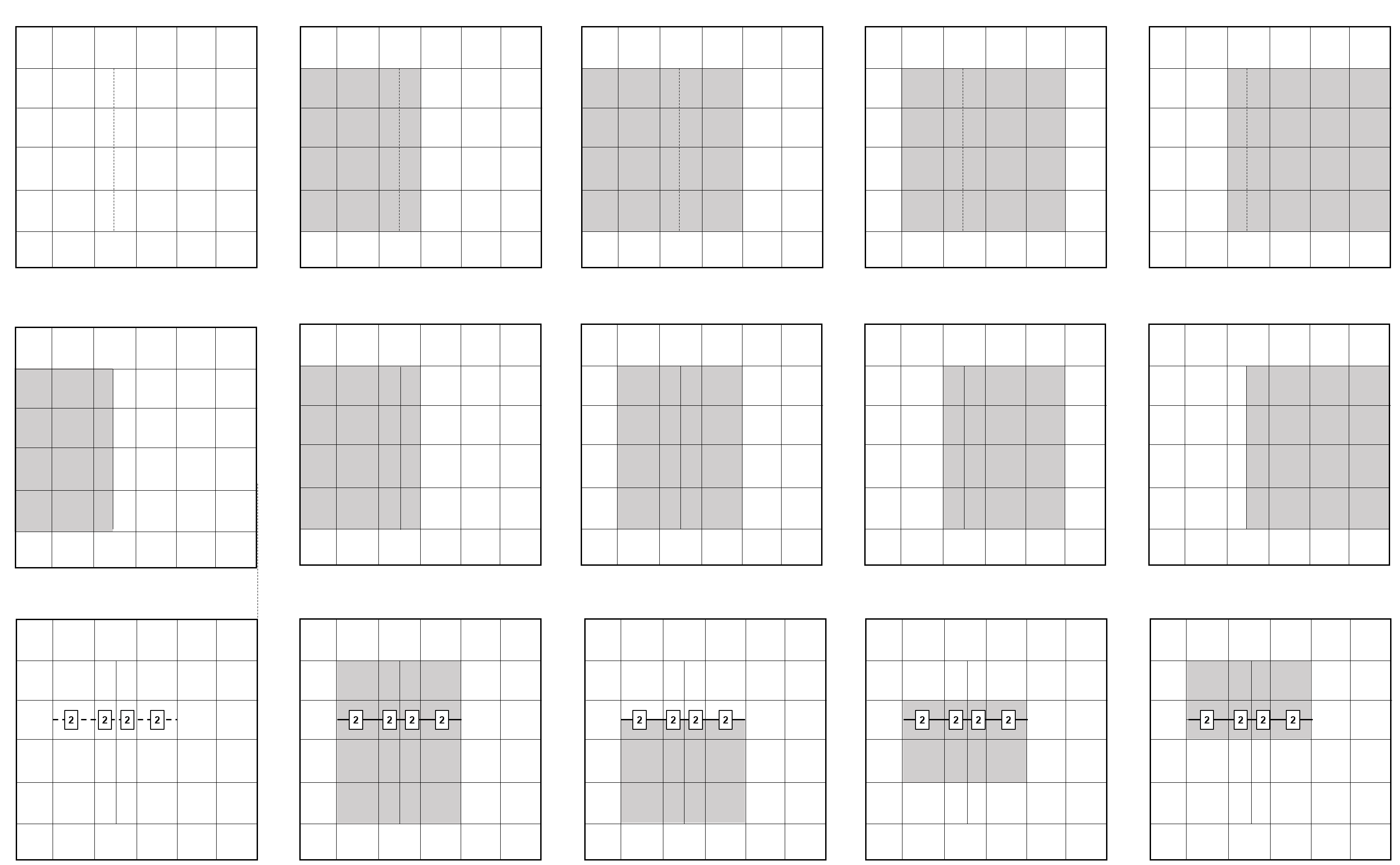}
\caption{On the upper left corner we have the mesh of a box-partition of a bi-degree (3,3) tensor product B-spline basis. We insert a meshrectangle that splits at least one B-spline over this mesh. The four meshes following to the right show the tensor product B-splines split by this refinement. In the middle row of meshes we show how the four tensor product B-splines are refined into five tensor product B-splines. In the left mesh on the last row we  insert one meshrectangle with multiplicity two. This only affects the B-spline highlighted in the second mesh on the last row. The next meshes show the resulting refined B-splines.
}\label{fig:LR_refinements}
\end{center}
\end{figure}

\subsection{LR B-splines}\label{sec-LR}

Locally Refined B-splines (LR B-splines) were introduced in 2013 \cite{Dokken-2013} and the use of LR B-splines in IGA presented in \cite{Johannessen-14}. The definition of LR B-splines is based on what is denoted box-partitions of a $d$-box in $\RR^d$, $d \ge 1$. The approach of LR B-splines is thus not restricted to the bivariate, or trivariate case. There is no restriction on the polynomial degrees to be used and odd and even degrees are treated in the same way. Knots can be inserted at arbitrary values, thus supporting non-uniform refinements. The refinement is performed by inserting a meshrectangle that splits at least one B-spline. In the bivariate case meshrectangles are knot-line segments, and in the univariate case they are knots. A meshrectangle can have multiplicity higher than 1, thus generalizing knot multiplicity of univariate B-splines. Rather than go deep into technical detail, we will illustrate the concept with examples of bi-degree (3,3) LR B-splines. For technical details, consult the cited publications.

In Fig. \ref{fig:B-spline_to_box-partition} we show how a the parameter plane and the knot values of a bi-degree (3,3) tensor product B-spline space is represented as a box partition with multiplicities assigned to the mesh rectangles (knot-line segments). The refinement of LR B-splines is illustrated in Fig. \ref{fig:LR_refinements}. The rules for LR B-spline refinement are as follows:
\begin{itemize}
\item Insert a meshrectangle that splits the support of at least one B-spline.
\item Find all B-splines that have a support that is split by the meshrectangle. Refine these B-splines until the collection of B-splines contains only minimal support B-splines. Note that existing meshrectangles might split B-splines resulting from the refinement. Such B-splines must be recursively subdivided until all B-splines have minimal support. 
\end{itemize}

As with T-splines there are refinement configurations that can result in a linearly dependent collection of tensor product B-splines. However, linear independence can be ensured if a hand-in-hand refinement condition is fulfilled throughout the refinement process. This means that the dimension of the spline space is checked against the increase in number of B-splines at each refinement step.

\begin{figure}[h]
\begin{center}
\includegraphics[width=10.0 cm]{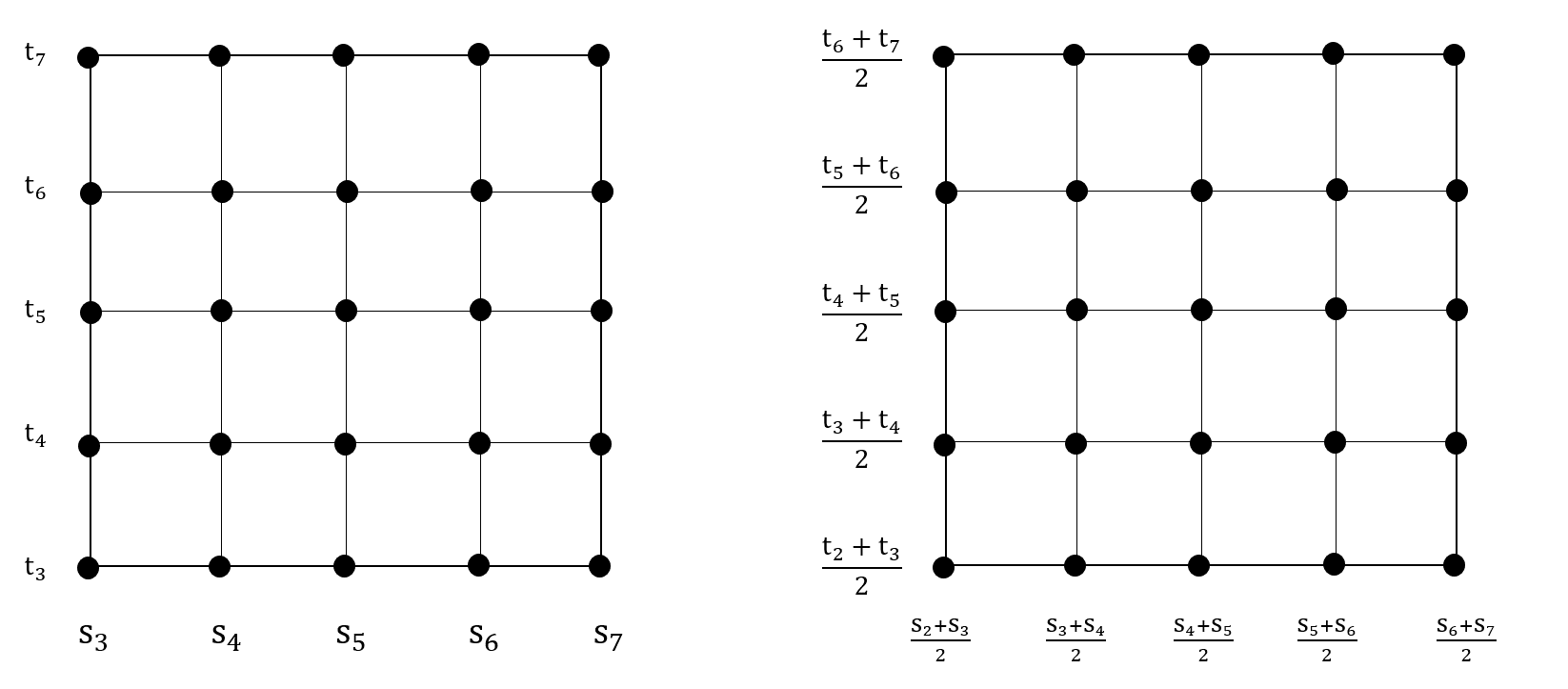}
\caption{To the {\it left} we show how knot values can be assigned to the control points of tensor product B-splines of bidgree $(3,3)$. To the {\it right} we show how knot values can be assigned to control points of B-splines of bidegree $(2,2)$ by averaging the middle knots thus creating the Greville point of the tensor product B-spline.
}\label{fig:LR-vertex_refinement}
\end{center}
\end{figure}

\section{Comparing Hierarchical B-splines, T-splines and LR B-splines}\label{sec-comparing}
In the previous sections we have explained the ideas behind (truncated) hierarchical B-splines, T-splines and LR B-splines. In this section we will try to explain how the different approaches relate. Can one method mimic the properties of the other approaches? In Section \ref{sec-specifying-the-refinement} we will investigate the strategies for specifying the refinement. This is followed in Section \ref{sec-degree-dim} by discussion of generalization in degree and dimension, and then in Section \ref{sec-spline-spaces} the difference of spline spaces will be addressed. Linear independence is important when doing isogeometric analysis and this is addressed in Section \ref{sec-linear-independence}.

\subsection{Refinement specification}\label{sec-specifying-the-refinement}
Hierarchical B-splines, T-splines and LR B-splines have different strategies for specifying the refinement to be performed, respectively specification of the region to be refined for the next hierarchical level, inserting new control points in the T-mesh, and inserting meshrectangles in the box-partition. A natural question to pose is if these ``user interfaces'' can be interchanged.

\begin{itemize}
\item {\bf Specify region of interest for next level of hierarchical refinement.} As both T-splines and LR B-splines have very high flexibility with respect to refinement it is feasible to impose a hierarchical B-spline type of refinement rule to both methods and specify refinement level-by-level. The tensor product B-splines to be refined can be selected by a specification of regions as for HB/THB. For LR B-splines this hierarchical B-splines approach has been followed using what we call structured refinement \cite{Johannessen-14}, where all elements of the selected B-splines are split in two in all parameter directions at each refinement level. For bi-degree (3,3) and lower such LR B-splines are linearly independent, but linear dependence has been observed in degrees $\geq4$.
\item {\bf Refinement by insertion of new control points in the control point mesh.} As hierarchical B-splines are intrinsically defined by the region of subdivision, this is not a feasible approach. However, the control point insertion of T-splines can be used for LR B-splines. In Fig. \ref{fig:LR-vertex_refinement} to the {\it left} we show that for LR B-splines of bidegree (3,3) the control points can be anchored to the middle knot pair of the corresponding tensor product B-spline similar to the anchoring of control points for T-splines. In Fig. \ref{fig:LR-vertex_refinement} to the {\it right} we show that the control vertices of a bidegree (2,2) tensor product B-spline/LR B-spline surface can be assigned a parameter value pair corresponding to the Greville point. Note that after splitting such Greville point based parameter value pairs must be updated. 
\item {\bf Refinement by inserting meshrectangles.} While LR B-splines enforce that all tensor product B-splines are minimal support after each refinement step, T-splines only enforce minimal support in the parameter direction of the refinement, and only on the B-splines split by the new control point. So, the LR B-spline approach can be given a very T-spline like behaviour if we relax the minimal support requirement to only be in the refinement direction, and ensure that refinements are between anchored control points. 
\end{itemize}

\subsection{Generalization in degree and dimension}\label{sec-degree-dim}
The approach of (truncated) hierarchical B-splines and LR B-splines are neither restricted in degree nor dimension. Both perform refinement in the parameter domain, and all degrees are treated in the same way. As T-splines perform refinement in the control mesh, more advanced navigation has to take place for identifying knot values as the number of dimensions increases. There are already examples of trivariate T-splines \cite{Solid-T-spline}. T-splines are explained using bi-degree (3,3), but the anchoring of tensor product B-splines to vertices works well for any odd degree. For even degrees the dual grid can be used.

\subsection{Differences of spline spaces for hierarchical type meshes}\label{sec-spline-spaces}
As the algorithms for finding which tensor product B-splines to use are different for truncated hierarchical B-splines, T-splines and LR B-splines the spline spaces generated will be different. However, when T-splines and LR B-splines mimic the knotline meshes of hierarchical B-splines, the behaviour is quite similar. In \cite{Johannessen-16} it is shown that on similar grids, LR B-splines and truncated hierarchical B-splines are fairly similar with respect to the condition numbers of the stiffness matrix for the examples considered. It must be expected that T-splines will behave in a similar way. However, the original formulation of hierarchical B-splines has a significantly higher condition number due to the lack of partition of unity of the resulting B-splines. In general, the number of tensor product B-splines will grow much faster for (truncated) hierarchical B-splines than LR B-splines and T-splines. For LR B-splines and T-splines a single vertex/knotline segment can be inserted, while for (truncated) hierarchical B-splines refinements dictate that a much higher number of knotline segments will be added. Consequently, for large problems it much be expected that LR B-splines and T-splines will outperform truncated hierarchical splines.

\subsection{Linear independence}\label{sec-linear-independence}
Truncated Hierarchical B-splines are always linearly independent. LR B-splines can be defined over most meshes of hierarchical B-splines (in some cases extra refinement is necessary). Linear independence is guaranteed for LR B-splines of bi-degree (3,3) or lower over such hierarchical meshes. The only example known for T-splines that is linearly dependent has multiplicity of knot values in the interior of the mesh. Consequently, we can expect that T-splines on meshes similar to meshes of hierarchical B-splines are linearly independent. For LR B-splines of higher bi-degree than (3,3) on hierarchical meshes a simple rule can be imposed on the mesh that guarantees linear independence: if the region to be refined can be split into two overlapping subregions whose intersection does not contain any of the B-splines to be refined then special care has to be taken to check if the region has to be extended to avoid linear dependence.

Due to the flexibility of LR B-splines and T-splines there are situations where the resulting collection of tensor product B-splines are linearly dependent. For bivariate LR B-splines it is always possible to use the hand-in-hand principle mentioned above to ensure that a refinement will produce a linearly independent set of B-splines. However, this depends on the availability of an appropriate dimension formula, which is currently only available in the bivariate case. Work is going on to find other subclasses of LR B-spline that ensure linear independence without using the hand-in-hand principle. For T-splines the subset of analysis suitable T-splines is always linearly independent.

\section{Trivariate spline extensions to ISO 10303} \label{Sec-tri-var-cad}
The development of B-spline technology in the 1970s and NURBS in the 1980s had a rapid uptake both in CAD-industry and in the 1990s by ISO 10303 STEP \cite{ISO 10303-28:2007}. Prior representations for freeform curves and sculptured surfaces were replaced by curves and tensor product surfaces based on B-spline and NURBS technology. As the subtractive and formative manufacturing technologies at that time were based on uniform material there was no need for representing the interior of an object. 

Before IGA, addressed in Section \ref{Sec-IGA}, was accepted as a technology with great industrial potential, trivariate B-spline representations for volume objects attracted little attention in STEP. The consequence was that the trivariate B-spline representation already in Part 42 \cite{ISO 10303-42:2016} was not described in Application Reference Model (ARM) Part 1801 {\it B-spline Geometry}, see Appendix \ref{App-Extended-B-spline-geometry} for details. The ARMs of STEP are written as a reference for those that develop converters. 

The EC Factories of the Future project TERRIFIC (2011-2014, Contract No. 284981) proposed several additions related to Locally Refined Splines to STEP Part 42. In addition, extensions of STEP AP 242 Edition 2 \cite{STEP-AP242ed2} were proposed. The extensions to Part 42 will be published in 2018. To support the Part 42 extension the preparation of a new ARM, currently denoted {\it Extended B-spline Geometry}, was started in 2018 in the EC Factories of the Future project CAxMan (2015-2018, Contract No. 680448). For details see Appendix \ref{App-Extended-B-spline-geometry}. To support trivariate spline representations, extensions are needed in other parts of STEP such as Part 50 and Part 52. These extensions are carried out by other parts of the STEP-community.

\section{Analysis-based design for AM} \label{Sec-ana-bas-des}
Additive manufacturing is moving from objects composed of a single material to objects composed of multiple discrete materials and graded materials. The additive manufacturing processes create objects with material properties that are significantly more anisotropic than traditional manufacturing technologies. Consequently, it will be beneficial to use extensive analysis already during the design stage to better support the additive manufacturing technology chosen. 

The specificities of the different additive processes have a strong influence on how an object should be designed. For example, for metal powder bed based additive technologies, overhangs of more than 45 degrees require support structures to be added. However, if the printing direction is known then the ``roof'' of an internal void or cavity could possibly be replaced by a drop shaped roof thus avoiding overstepping the 45-degree restriction. However, the consequences of such design changes targeting manufacturability should be analysed. 

It is well known that meshing from CAD is work intensive. The DART study of Sandia Labs quantified the increasing challenge of mitigating the problems with non-watertight CAD-models when creating suitable meshes for FEA \cite{Dart:2005}. Time spent fixing such issues was reported to have increased from 73\% in 2005 and numbers as high as 90\% have been mentioned in 2017. The large CAD vendors invest much effort to make B-rep CAD and FEA seamlessly interoperable in their market offerings. 

To go beyond what the large vendors offer, it seems that a suitable approach is to build analysis-based design for AM on a combination of trivariate CAD-models and IGA. Originally IGA was based on block-structured trivariate spline models, an approach like CAD-models before trimming of B-spline surfaces was introduced. To augment the trivariate spline representation it is natural to trim the trivariate spline models by surfaces, thus generalizing trimming of spline surfaces by loops of edges to trimming of spline volumes by B-rep shells. This is the same line of ideas as the embedded methods of cut-FEM and the Finite Cell Method (FCM). However, to trim the trivariate spline model with a CAD B-rep model that is not watertight is as far as we know new. In Section \ref{Sec-trimmed-iga} we address the work performed in the CAxMan project on trimmed trivariate spline models.

\subsection{Isogeometric Analysis} \label{Sec-IGA}In 2005 Tom Hughes \cite {Hughes:05} introduced Isogeometric Analysis (IGA). IGA builds on the ideas of FEA. In FEA the shape functions of an element are local to the element. In IGA the shape functions are replaced by B-splines, that cross element boundaries. Thus, the same B-spline can be used as shape function in more than one element. The transition between elements will replicate the continuity of the B-splines used. As B-splines and Non-Uniform Rational B-splines can represent all CAD-shapes, IGA can in principle have an accurate representation of all shapes used in CAD-design \cite{Skytt:2016}. B-splines can have any polynomial degree, consequently IGA supports the use of higher degree element representations in simulation. IGA has been demonstrated to be more accurate than traditional FEA. However, making a trivariate block-structured spline model suited for IGA from a CAD-model meets the similar challenges to making a block-structured FEA-mesh from CAD.

Advantages and challenges of block-structured IGA (no trimming):
\begin{itemize}
\item Easy to handle in analysis.
\item An approximation step is normally required in the block structuring process.
\item Difficult to create IGA block structured models for complex objects:
\begin{itemize}
\item For surfaces models there might occur vertices with valence different from four, e.g., less than four or more than four faces connect at a vertex. 
\item For surfaces models there might occur singular edges.
\item For volume models there might occur edges with valence different from four, e.g., less than four or more than four topological volumes connect along an edge.
\item For volume models there might occur singular edges and faces.
\end{itemize}
\end{itemize}

\begin{figure}[t]
\begin{center}
\includegraphics[width=13.5 cm]{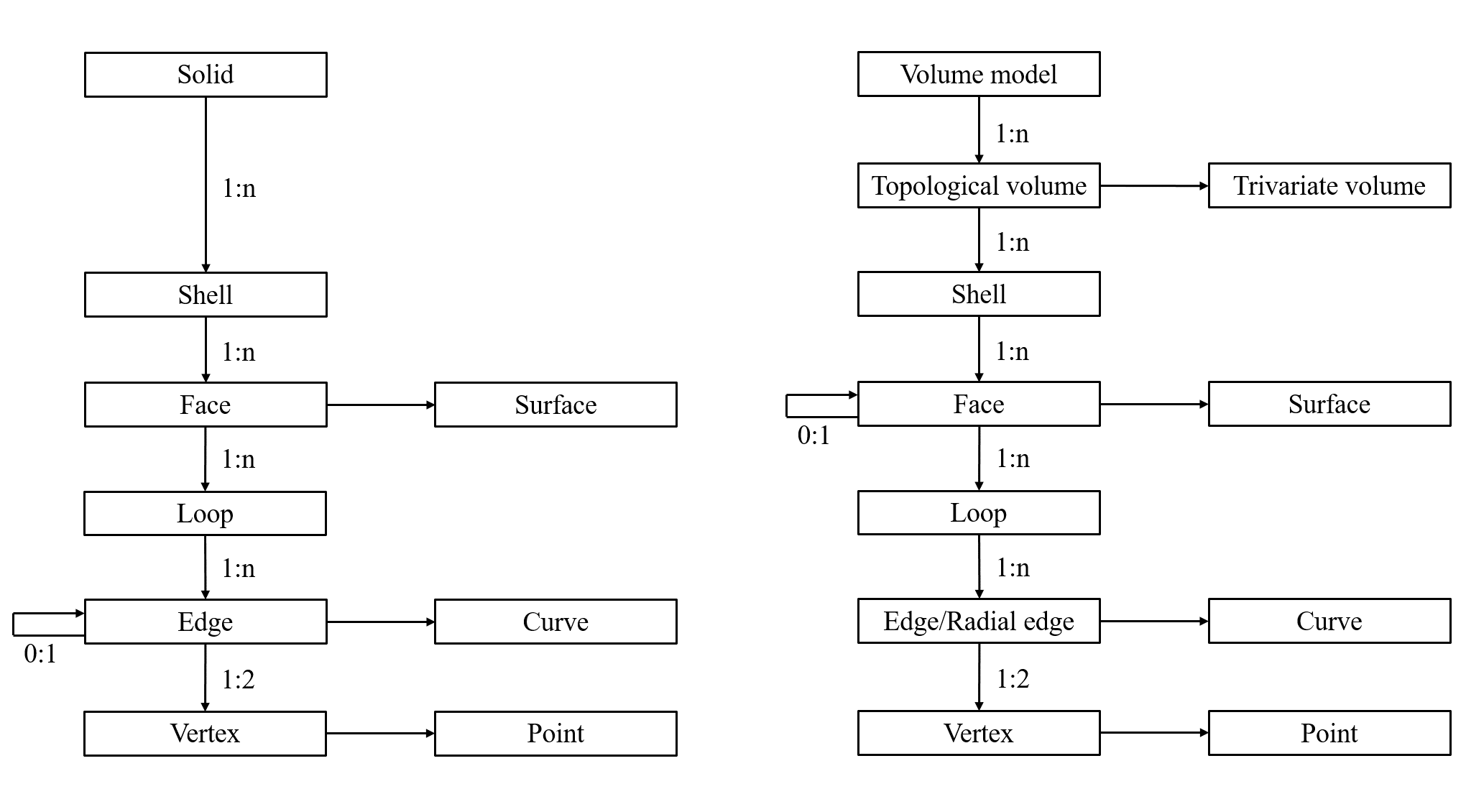}
\caption{To the {\it left} a CAD type B-rep topology data structure, and to the {\it right} a data structure addressing the needs of a block-structured model for IGA where the blocks may be trimmed volumes. In a B-rep model a face is unique, and the topological outline(s) of the face is described by loop(s) of edges. When two faces are adjacent the connection is described through relations between coincident edges of the two faces (indicated by the $0:1$ on the relation between edges in the left picture). In the trimmed trivariate block-structure for IGA, two topological volumes that are adjacent are connected through coincident faces through a relation (indicated by the $0:1$ on the relation between faces in the right picture). While an edge in a B-rep model can only be related to two faces, an edge in the trimmed trivariate block-structure for IGA can be coincident to edges of multiple faces. 
}\label{fig:Brep_and_trivariate}
\end{center}
\end{figure}

\subsection{Trimmed trivariate CAD, IGA and AM} \label{Sec-trimmed-iga}

In Fig. \ref{fig:Brep_and_trivariate} we show to the left a typical topological structure for B-rep CAD, and to the right a proposed topological structure for trimmed trivariate CAD \cite{Skytt:2016}. In the latter, two new entities are introduced: the topological volume and the trivariate volume that contains the mathematical representation. This approach allows the B-rep CAD-model to trim the trivariate volumes. Although this in principle seems simple, the difference between what is regarded as a watertight CAD model and as a watertight FEA or IGA-model remains.
\begin{itemize}
\item A CAD model is regarded as watertight if the topology representation is correct and the gaps between adjacent faces are within user defined tolerances. 
\item A FEA model is watertight when the elements match exactly, no gaps are allowed between elements. 
\item A volumetric IGA model is watertight when spline volumes match exactly, no gaps are allowed between adjacent volumes. 

\end{itemize}

So, to trim a block-structured trivariate IGA model with B-rep CAD-models poses algorithmic challenges. If the original B-rep model includes gaps and these live on through subsequent uses of the model, algorithms must always take the challenges of gaps into consideration. Alternatively, the model must be repaired to remove gaps, and thus simplify subsequent uses.
In the CAxMan project, the focus is on trimmed trivariate block-structured models for IGA and AM. 

\begin{figure}[t]
\begin{center}
\includegraphics[width= 13.5cm] {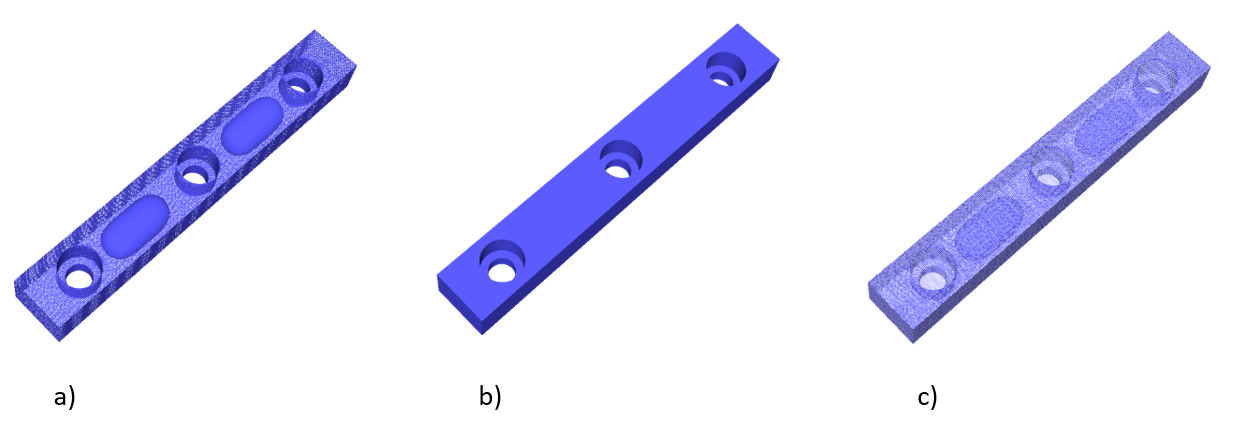} 
\caption{In a) a transparent view showing the voids. In b) the outer shell of a CAD-model.  In c) the trimmed trivariate spline representation, where the white shadow follows the main shape of the volume.}
\label{fig:CAD-and-trimmed-trivariate}
\end{center}
\end{figure}

\begin{figure}[h]
\begin{center}
\includegraphics[width= 13.5cm] {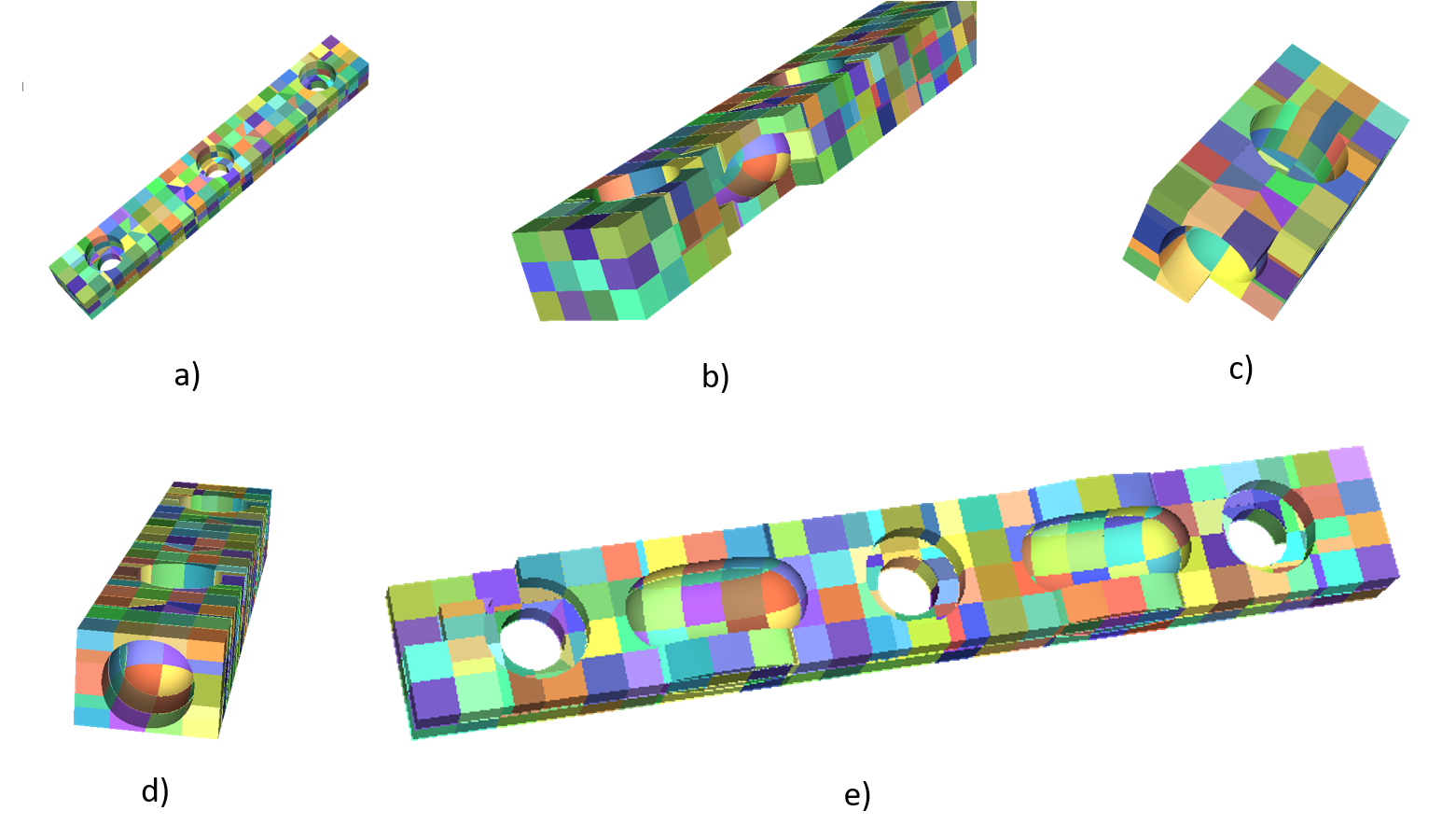} 
\caption{In image a) the detrimming and reparameterization to trivariate hexahedra shapes of the trimmed trivariate model. In b) some spline hexahedra are removed to show the reparameterization around the void, and in c) and d) a void seen from another direction with other spline hexahedra removed. In e) a section of the detrimming is seen from the top.} \label{fig:Detrimming_illustrations}

\end{center}
\end{figure}
\begin{figure}[h]
\begin{center}
\includegraphics[width=13.5cm]{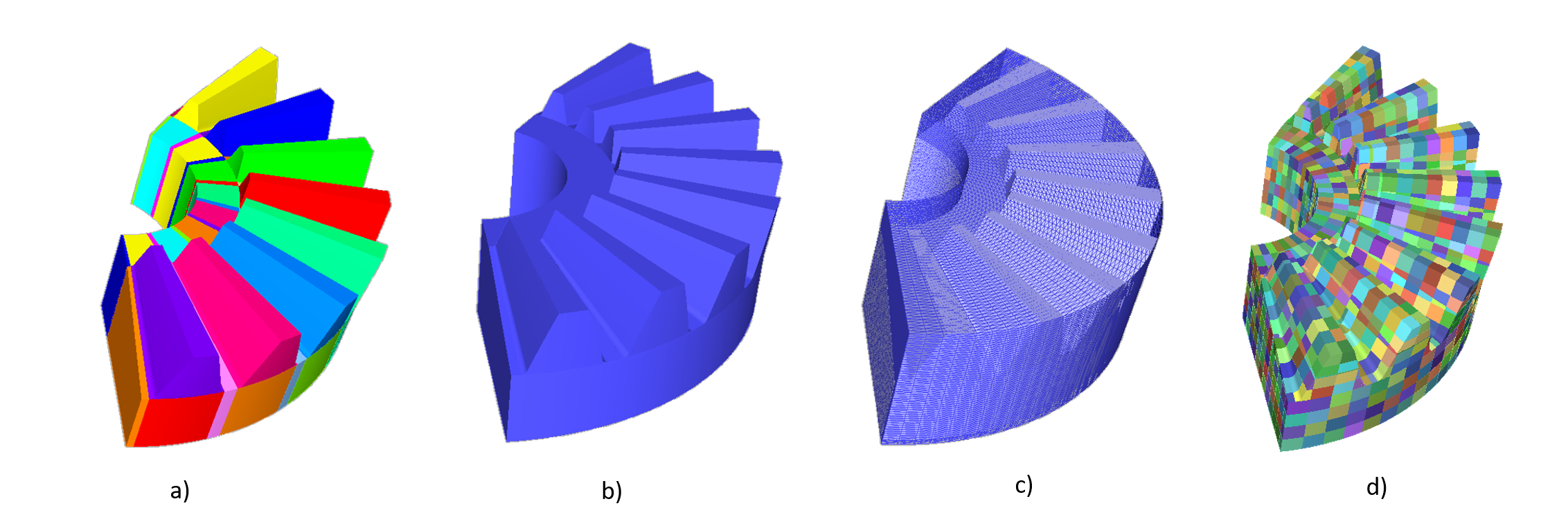}
\caption{In a) a block-structured model of a gear. The same gear represented as a trimmed spline volume in b), then in c) the trimmed volume with the trimming surfaces, and d) is the detrimming prepared for quadrature. Note that the shape of the trivariate trimmed spline follows the rotational shape of the gear and thus ensures rotational symmetry. The example is courtesy of the CAxMan partner Stam S.r.l.. \label{fig:trimmed_examples}}
\end{center}
\end{figure}

\section{Detrimming of trimmed trivariate spline models}\label{Sec:quadrature}
Quadrature of a trivariate spline volume that is parameterized over a hexahedron can be decomposed into univariate quadrature in each of the parameter directions. Consequently, to be able to perform numerical integration over the trimmed trivariate spline model, it makes sense to split the model into a collection of sub-volumes each parameterized over a hexahedron. In CAxMan we have implemented a first version of such a {\it detrimming algorithm}. Fig. \ref{fig:CAD-and-trimmed-trivariate}.a) illustrates the conversion of a B-rep CAD-model with voids to a trimmed trivariate spline model. This is then detrimmed for quadrature as illustrated in Fig. \ref {fig:Detrimming_illustrations}a--e.

Fig. \ref{fig:trimmed_examples} illustrates for a part of a gear, a pure block-structured model as well as a trimmed block-structured model of the part of the gear, and the detrimming of this.

\section{Open challenges}\label{Sec:challenges}
To represent accurately the fine detail of objects with large interior lattice structures, B-rep CAD is not feasible as the number of geometric elements will be much larger than what was foreseen when this format was established. In addition to the bulk of the number of geometric elements, geometric tolerances will also pose a challenge. In general B-rep based CAD-systems allow gaps between adjacent surfaces controlled by user defined tolerances. One solution to this is to set finer tolerances. However, then CAD models that are valid with their original tolerance settings will be invalid with the finer tolerances.

When designing lattice structures, the size of a gap between adjacent faces and the thickness of the volumetric elements of the lattice structures risk not having a proper separation. This potentially creates failure of the CAD-system. This might be one of the reasons that lattice structures are today frequently created after the conversion to STL, with the consequence that the CAD-model of the object manufactured does not accurately mirror the object manufactured.

The use of FEA or trivariate spline based CAD-model representation will solve the tolerance problem \cite{Weeger:16} \cite{Ezair:17} \cite{Massarwi:16}. However, the bulk of the representation will grow. Automatic generation of trimmed trivariate CAD-models is still a technology under development, so there are still significant developments necessary before it is as mature as B-rep based CAD. More innovative approaches to lattice structure representations must be pursued in the future. At SIM-AM in Munich in October 2017 many different promising approaches for modelling of lattice structures were presented. However, the challenge is to ensure that these can be represented in a CAD-based digital twin that reflects the object manufactured by additive technology. Our feeling is that a lot seems to remain before we have a generic technology for representation and modelling of lattice structures for AM.

\section{Acknowledgement}
This project has received funding from the European Union's Horizon 2020 research and innovation programme under grant agreement No 680448.





\begin{appendices}
\chapter{Extensions to ISO 103030 (STEP)}
\section{Application Reference Model 1801: B-spline Geometry}\label{App-B-spline-geometry}
ARM 1801 includes the following well know representations of geometry using B-splines:
\begin{itemize}
\item {B-spline curve.}
\item {B-spline surface.}
\item {Rational B-spline curve.}
\item {Rational B-spline surface.}
\item {Surface with explicit knots.}
\item {Surface with implicit knots.}
\end{itemize}

\section{Proposal for new ISO 10303 Application Reference Model: Extended B-spline Geometry}\label{App-Extended-B-spline-geometry}
The ARM to be proposed on {\it Extended B-spline Geometry} is planned to include:
\begin{itemize}
\item {\bf B-spline volume and B-spline volume with knots.} Although these have been part of Part 42 for a long time, an ARM description has not been made until now, probably due to little interest in its use before IGA emerged. This is a parametric volume represented by a trivariate tensor product spline basis.
\item {\bf Status of linear independence} is important for LR B-splines and T-splines. The values this variable can have are {\it Independent}, {\it Not independent}, and {\it Not tested}. 
\item {\bf List of types of Locally Refined Splines.} The list of values for this variable currently includes {\it Analysis Suitable T-spline}, {\it Hierarchical B-spline}, {\it LR B-spline}, {\it Semi-Standard T-spline} and {\it Standard T-spline}. The Truncated Hierarchical B-splines are currently not included but can potentially be represented by expanding the truncation and assigning the vertex values to all vertex values used, including the occurrences used in the truncation. This will result in multiple occurrences of the same B-spline as all its occurrences will be explicitly represented. 
\item {\bf Local B-spline.} As both Hierarchical B-splines, T-splines and LR B-splines are based on collections of (tensor product) B-splines with control points and scaling factors, each individual (tensor product) B-spline has to be represented as an entity. Each (tensor product) B-spline entity points for each dimension to a list of knot values with a parallel list of multiplicities.
\item {\bf Locally refined spline curve and rational locally refined spline curve.} This entity allows us to extracted constant parameter line spline curve from a locally refine spline surface of volume represented by B-splines from the locally refine spline surface of volume. In general, this description will not be a minimal support set of B-splines, and will be represented by a collection of B-splines that risk being linearly dependent. However, these spline curves can be exactly converted to a (Rational) B-spline curve that is using a linearly independent collection of B-splines when necessary. 

For many uses it can often be wise to convert the curve to a minimal support B-spline basis, i.e., a B-spline basis using a normal knot vector.
\item {\bf Locally refined spline surface and rational locally refined spline surface} can be of types listed in {\it List of types of Locally Refined Splines}. The surface is represented by a collection of bivariate tensor product B-splines with scaling factors and control point values.

\item {\bf Locally refined spline volume and rational locally refined spline volume } can be of types in {\it List of types of Locally Refined Splines}. The volume is represented by a collection of trivariate tensor product B-splines with scaling factors and control point values.
\end{itemize}

\end{appendices}

\end{document}